\documentclass{amsart}

\copyrightinfo{2003}{American Mathematical Society}

\newtheorem{theorem}{Theorem}[section]
\newtheorem{lemma}[theorem]{Lemma}
\newtheorem{proposition}[theorem]{Proposition}
\newtheorem{corollary}[theorem]{Corollary}

\theoremstyle{definition}

\theoremstyle{remark}

\numberwithin{equation}{section}

 \usepackage{graphicx}

\usepackage{amssymb}

\newcommand{\Nset}{\mathbb{N}}
\newcommand{\Rset}{\mathbb{R}}
\newcommand{\de}{\,\mathrm{d}}
\newcommand{\e}{\mathrm{e}}

\begin{document}

\title[On The Boundedness In $H^{1/4}$ Of The Maximal Square Function]
{On The Boundedness In $H^{1/4}$ Of The Maximal Square Function
Associated With The Schr\"odinger Equation}

\author{Giacomo Gigante}
\address{Dipartimento di Ingegneria Gestionale e dell'Informazione,
 Universit\`a di Bergamo,
Viale Marconi 5, 24044 Dalmine (BG), Italy.}
\email{gigante.giacomo@unibg.it}

\author{Fernando Soria}
\address{Departamento de Matem\'aticas, Facultad de Ciencias, C-XV,
Universidad Aut\'onoma de Madrid, 28049 Madrid, Spain.}
\email{fernando.soria@uam.es}
\thanks{This research was partially supported by the European
Commission via the Harmonic Analysis Network ``HARP'', and by grant BFM2001-0189.}

\subjclass[2000]{42B15, 42B25}
\keywords{Maximal square function, linear Schr\"odinger equation, Sobolev spaces.}

\date{July 9, 2003}

\begin{abstract}
A long standing conjecture for the linear Schr\"odinger equation states that
$1/4$ of derivative in $L^2$,
in the sense of Sobolev spaces, suffices in any dimension for the solution to that
equation to converge almost everywhere to the initial datum as the time goes to 0.
This is only known to be true in dimension 1 by work of Carleson. In this paper
we show that the conjecture is true on spherical averages. To be more precise,
we prove the $L^2$ boundedness of the associated maximal square function on the
Sobolev class $H^{1/4}(\Rset^n)$ in any dimension $n$.
\end{abstract}

\maketitle

\begin{section} {Introduction}
For $\alpha\in \Rset$, we denote by $H^\alpha(\Rset^n)$ the Sobolev space
\[
H^\alpha(\Rset^n)=\left\{f\in{\mathcal S}^\prime(\Rset^n): \|f\|_{H^\alpha}=
\left(\int\left|\widehat f (\xi)\right|^2(1+|\xi|^2)^\alpha\de\xi\right)^{1/2}
<\infty \right\}.
\]
We will also consider the homogeneous Sobolev space $\dot{H}^\alpha(\Rset^n)$
defined by
\[
\dot{H}^\alpha(\Rset^n)=\left\{f\in{\mathcal S}^\prime(\Rset^n): \|f\|_{\dot{H}^\alpha}=
\left(\int\left|\widehat f (\xi)\right|^2|\xi|^{2\alpha}\de\xi\right)^{1/2}
<\infty \right\}.
\]
Let $f$ be in the Schwartz class ${\mathcal S}(\Rset^n)$, and define
\[
S_tf(x)=u(x,\,t)=\int_{\Rset^n}\widehat f(\xi)\e^{-2\pi \mathrm{i}|\xi|^2t}
\e^{2\pi \mathrm{i} \xi\cdot x}\de\xi.
\]
Then $u$ is the solution to the linear Schr\"odinger
equation with initial datum $f$, that is,
\[
\left\{
    \begin{array}{ll}
        \frac\partial{\partial t}u(x,t)=\frac {\mathrm{i}}{2\pi}\Delta_xu(x,t)
            & \mbox{in $\Rset_+^{n+1}$}\\
        u(x,0)=f(x) & \mbox{in $\Rset^n$.}
    \end{array}\right.
\]

There is a fundamental question in this setting and is that of determining the minimal
smoothness on the initial value function $f$, needed for the almost
everywhere convergence
\begin{eqnarray}
\label{convergenza}
\lim_{t\to 0^+}u(x,\,t)=f(x),\qquad\mbox{a.e.}
\end{eqnarray}
This smoothness is measured in terms of the Sobolev space $H^\alpha$ which
the function $f$ belongs to. In 1979, Carleson proved in \cite{carleson} that
the $\mbox{a.e.}$ convergence (\ref{convergenza}) holds for any $f\in\dot{H}^{1/4}$
in dimension $n=1$. Dahlberg and Kenig \cite{kenig} extended this result to functions
in $\dot{H}^{n/4}(\Rset^n)$ for any $n$ and showed that there are counterexamples
if the regularity is less than $1/4$. It is conjectured that $\alpha=1/4$ suffices
for this problem in any dimension $n$. Sj\"olin and Vega proved independently
in \cite{sjolin}, \cite{vega}
that $\alpha$ greater than $1/2$ implies the convergence (\ref{convergenza})
in any dimension $n$ (previous results,
for $\alpha>1$, were obtained in \cite{carbery}, \cite{cowling}), while Prestini
\cite{prestini} proved the conjecture for radial functions.
The case $n=2$ has been intensively studied during the last years and is the only
one (apart from $n=1$) where there are positive results for ($\ref{convergenza}$)
with smoothness
$\alpha<1/2$ (see \cite{vargas}, \cite {tao}, \cite{taovargas}, and the references there).

As usual, problems related to the $\mbox{a.e.}$
convergence are intimately connected to the boundedness of the associated
maximal function. In our case, this maximal function is given by
$S^\ast f(x)=\sup_{t>0}\left|S_t f(x)\right|$, $x\in\Rset^n.$
For example, the $\mbox{a.e.}$ convergence (\ref{convergenza}) for all
functions $f\in H^\alpha$ follows from the a priori maximal estimate
\begin{eqnarray}
\label{limitatezza}
\left(\int_{|x|\le 1}|S^\ast f(x)|^p\de x\right)^{1/p}\le C\|f\|_{H^\alpha},\qquad
f\in{\mathcal S}(\Rset^n).
\end{eqnarray}
In fact, all the known cases about convergence mentioned above are obtained
via this maximal inequality for different values of $p\in[1,\,2]$.

In this paper, we investigate whether inequality (\ref{limitatezza}) holds
if we replace $S^\ast$ by a spherical average operator; namely we look at the
maximal square function
\[
Q^\ast f(x)=
\sup_{t>0}\left(\frac 1{\sigma(S^{n-1})}\int_{S^{n-1}}\left|S_tf(|x|\omega)\right|^2
\de\sigma(\omega)\right)^{1/2}.
\]
Clearly, one has the inequality $\int_{|x|\le 1}
|Q^\ast f(x)|^2\de x\le\int_{|x|\le 1}|S^\ast f(x)|^2\de x $, and therefore
the boundedness of $S^\ast$ would imply a corresponding inequality for $Q^\ast$.
The known counterexamples show that the smoothness condition $\alpha\ge 1/4$ is
still necessary for the boundedness of this operator.
The main result of this paper says that $\alpha=1/4$
is also sufficient for the boundedness of $Q^\ast$.
\begin{theorem}
\label{main}
The operator $Q^\ast$ is bounded from $\dot{H}^{1/4}(\Rset^n)$ into $L^2(\{|x|\le 1\})$
in any dimension $n$; in fact, there is a positive
constant $C$, independent of the dimension, such that
\begin{eqnarray}
\label{stima}
\left(\int_{|x|\le 1}|Q^\ast f(x)|^2\de x\right)^{1/2}\le C\|f\|_{\dot{H}^{1/4}},\qquad
\forall \,f\in{\mathcal S}(\Rset^n).
\end{eqnarray}
\end{theorem}

In particular, (\ref{stima}) gives us that $1/4$ of smoothness suffices for the
$\mbox{a.e.}$ convergence with respect to quadratic spherical means. The precise
statement is contained in the following corollary.

\begin{corollary}
\label{corol}
If $f\in H^{1/4}(\Rset^n)$, then, for every $x_0\in\Rset^n$ we have
\[
\lim_{t\to 0^+}\int_{S^{n-1}}
|S_tf(x_0+r\omega)-f(x_0+r\omega)|^2\de\sigma(\omega)=0,\quad\mbox{a.e.
$r$.}
\]
\end{corollary}

\noindent{\bf Proof.} The proof is standard. By translation invariance, we
may assume without loss of generality that $x_0=0$. It is easy to
see that, if $g\in {\mathcal S}(\Rset^n)$, then $S_tg\to g$ as $t\to
0^+$, uniformly in $\Rset^n$. Given $f\in H^{1/4}(\Rset^n)$ we take a
sequence $\{g_k\}_{k=1}^\infty\subseteq{\mathcal S}(\Rset^n)$ such
that $g_k\to f$ in $H^{1/4}(\Rset^n)$. Denote by $\mu$ the Borel
measure $\de\mu(r)=r^{n-1}\de r$. Let $\lambda>0$,
$B^n=\{x\in\Rset^n\colon|x|\le 1\}$ and define
\[
A_\lambda=\left\{0<r<1 \colon \limsup _{t\to 0^+} \int
_{S^{n-1}}\left|S_tf(r\omega)-f(r\omega)\right|^2
\de\sigma(\omega)>\lambda\right\}.
\]
Then, for any positive integer $k$
\begin{eqnarray*}
\mu(A_\lambda) &\le&  \mu\left(\left\{0<r<1 \colon \limsup _{t\to
0^+} \int _{S^{n-1}}\left|S_t(f-g_k)(r\omega)\right|^2
\de\sigma(\omega)>\frac\lambda 2\right\}\right)\\
&+& \mu\left(\left\{0<r<1 \colon \int
_{S^{n-1}}\left|g_k(r\omega)-f(r\omega)\right|^2
\de\sigma(\omega)>\frac\lambda 2\right\}\right).
\end{eqnarray*}
Now, Chebyshev's inequality and Theorem \ref{main} imply that
\begin{eqnarray*}
\mu(A_\lambda) &\le&  \frac C\lambda\int _{B^n} \sup _{t>0}\int
_{S^{n-1}}\left|S_t(f-g_k)(|x|\omega)\right|^2
\de\sigma(\omega)\de x\\
&+& \frac C\lambda \|f-g_k\|^2_{2}< \frac C\lambda
\|f-g_k\|^2_{H^{1/4}}, \quad \forall\,\,k,
\end{eqnarray*}
and, therefore, $\mu(A_\lambda)=0$.
\qed

Before we proceed with the proof of Theorem \ref{main}, let us
first make a reformulation of our problem and some additional comments.
Observe that if $\{{\mathcal Y}_k\}$ is an orthonormal basis of spherical
harmonics in $L^2(S^{n-1})$, and
$\hat f(\xi)\sim\sum_k f_k(|\xi|){\mathcal Y}_k(\xi/|\xi|)$
denotes the corresponding expansion of $\hat f$ with respect to this basis, then
\[
Q^\ast f(x)=\sup_{t>0}\left(\frac{2\pi}{\sigma\left(S^{n-1}\right)}\sum_k
\frac 1{|x|^{n-1}}\left|Q^t_{\nu(k)}\left(f_k(s)s^{(n-1)/2}\right)(|x|)
\right|^2\right)^{1/2},
\]
where
\[
Q^t_\nu g(r)=\int_0^\infty \e^{-2\pi\mathrm{i}ts^2}\widetilde J_\nu(2\pi rs)g(s)\de s
\]
and $\nu(k)=(n-2)/2+\mbox{degree}({\mathcal Y}_k).$
Here, $J_\nu$ denotes the Bessel function of order $\nu$ and $\widetilde J_\nu(t)
=\sqrt tJ_\nu(t)$ for $t\ge 0$.
Using that the norm in $\dot{H}^{1/4}$ of $f$ with respect to the above
expansion is given by $\|f\|_{\dot{H}^{1/4}}=\sum_k\int_0^\infty|f_k(r)|^2r^{1/2}
r^{n-1}\de r$,
and ``cancelling out the $\sum$ signs'', the inequality
\[
\int_{|x|\le 1}|Q^\ast f(x)|^2\de x\le C\|f\|^2_{\dot{H}^{1/4}}
\]
is equivalent to the estimate
\[
\int_0^1\sup_{t>0}|Q^t_\nu g(r)|^2\de r\le C\int|g(r)|^2r^{1/2}\de r,
\]
uniformly in the index $\nu$ too.

We can now follow
Carleson's approach (see \cite{carleson}, \cite{kenig}).
First we linearize our maximal operator, by making $t$ into
a function of $r$, $t(r)$. Next we may assume that $g$ is supported on a
fixed interval $I$ (as long as the final constant $C$ is independent of
$I$). ``Moving'' the smoothness to the other side (that is, redefining
$g(r)r^{1/4}$ as $g$ again), we consider instead the linear operator
\[
T_\nu g(r)=\int_I\frac{\e^{\mathrm{i}s^2t(r)}\widetilde J_{\nu}(rs)}
{s^{1/4}}g(s)\de s.
\]
Then what we have to show is
\begin{equation}
\label{estim}
\int_0^1|T_\nu g(r)|^2\de r\le C\int_I|g(s)|^2\de s,
\end{equation}
with $C$ independent of $g\in L^2(I)$, of the interval $I$, of the measurable function
$t(r)$ and of $\nu\in\Nset/2$.

We want to point out that Theorem~\ref{main} gives as a
consequence the boundedness of the maximal Schr\"odinger operator
$S^\ast$ on radial functions in $\Rset^n$, with constant
independent of $n$. A close look at the above arguments will
convince us that both, Theorem~\ref{main} and this dimension-free
estimate are, in fact, equivalent.

Let us bring here a related result obtained by the authors. In \cite{GS2} it was proved
that the uniform estimate
\[
\int_I \e^{\mathrm{i}as^2}J_\nu(s) \frac{\de s}{s^\beta}=O(1),
\]
independent of $\nu\in\Nset/2$, the interval $I$ and $a\in\Rset$,
holds (for $\beta<1$) if and only if $\beta\ge 1/6$.
This expression appears in a natural way as the leading term (using the product formula
for Bessel functions) of the expansion of the kernel associated to $T_\nu T_\nu^\ast$
but replacing the ``smoothness'' $1/4$ by the generic smoothness $\alpha$ with
$2\alpha-1/2=\beta$. This could be interpreted as an indication that the uniform estimate
of the operators $Q_\nu$ by this method would only be possible on the class
$\dot{H}^{1/3}$ ($\alpha=1/3$ corresponds to the case $\beta=1/6$). Our theorem here
shows that an additional cancellation of the rest of terms in the expansion of the
kernel is possible so that, as Theorem \ref{main} says, the result holds indeed on
$\dot{H}^{1/4}$.

Continuing with the reduction of our problem, let us point out that by using
a $TT^\ast$ argument and the well known expansion
\[
\widetilde J_\nu(r)=\sqrt{\frac 2\pi}\cos\left(r-\frac{\pi\nu}2-\frac\pi 4\right)+
O_\nu\left(\frac 1 r\right)\qquad\mbox{as $r\to\infty$,}
\]
it is not difficult to obtain (\ref{estim}) but with a constant which would depend
on $\nu$ (see also \cite{prestini}). Thus we only need to check that the
constant $C$ is uniformly bounded as $\nu$ tends to infinity.

The following lemma, due to J. A. Barcel\'o (\cite{BarTes}, \cite{BarRuVeg}),
describes the oscillation and the asymptotics of the Bessel function for large values,
with the precise dependence of the remainder term with respect to the order of the
function.
\begin{lemma}
\label{barcelo}
There is a universal constant $C>0$  such that for all
$\nu>1/2$ and for all $r>\nu+\nu^{1/3}$ we have
    \[
    J_\nu(r)=\sqrt{\frac{2}{\pi}}\frac{\cos \theta(r)}{(r^2-\nu^2)^{1/4}}
    +h_{\nu}(r),
    \]
where
    \[
    \theta(r)=(r^2-\nu^2)^{1/2}-\nu\arccos\frac{\nu}{r}-\frac{\pi}{4},
    \]
and
    \[
    |h_\nu(r)|\le\left\{
    \begin{array}{ll}
        C\left(\frac{\nu^2}{(r^2-\nu^2)^{7/4}}+\frac{1}{r}\right)
            & \mbox{ if $\nu+\nu^{1/3}\le r\le 2\nu$} \\
        \frac{C}{r} & \mbox{ if $r\ge 2\nu$.}
    \end{array}
    \right.
    \]
\end{lemma}
In order to simplify the notation, let us define for $r>\nu+\nu^{1/3}$ the functions
\begin{eqnarray*}
J_\nu^B(r)&=&\sqrt{\frac{2}{\pi}}\frac{\cos \theta(r)}{(r^2-\nu^2)^{1/4}},\\
\widetilde J_\nu^B(r)&=&\sqrt{r}J_\nu^B(r),\\
\widetilde h_\nu(r)&=&\sqrt{r}h_\nu(r).
\end{eqnarray*}
Thus, we can write $T_\nu$ as the sum of the following operators
\begin{eqnarray*}
T_\nu^1 g(r)&=&\int_I \e^{\mathrm{i}t(r)s^2}\widetilde J_\nu(rs)
\chi_{[0,\,\nu]}(rs)s^{-1/4}g(s)\de s,\\
T_\nu^2 g(r)&=&\int_I \e^{\mathrm{i}t(r)s^2}\widetilde J_\nu(rs)
\chi_{[\nu,\,\nu+\nu^{2/3}]}
(rs)s^{-1/4}g(s)\de s,\\
T_\nu^3 g(r)&=&\int_I \e^{\mathrm{i}t(r)s^2}\tilde h_\nu(rs)
\chi_{[\nu+\nu^{2/3},\,2\nu]}(rs)
s^{-1/4}g(s)\de s,\\
T_\nu^4 g(r)&=&\int_I \e^{\mathrm{i}t(r)s^2}\widetilde J_\nu^B(rs)
\chi_{[\nu+\nu^{2/3},\,2\nu]}(rs)
s^{-1/4}g(s)\de s,\\
T_\nu^5 g(r)&=&\int_I \e^{\mathrm{i}t(r)s^2}\tilde h_\nu(rs)
\chi_{[2\nu,\,\infty)}(rs)
s^{-1/4}g(s)\de s,\\
T_\nu^6 g(r)&=&\int_I \e^{\mathrm{i}t(r)s^2}\widetilde J_\nu^B(rs)
\chi_{[2\nu,\,\infty)}(rs)
s^{-1/4}g(s)\de s.\\
\end{eqnarray*}
The desired boundedness will now follow from the boundedness of the above operators.
This will be proved in sections 2 through 6, but first we would like to recall Van der
Corput's lemma.
\begin{lemma}[Van der Corput]
\label{vander}
Let $\phi$ be a smooth real valued function defined on an
interval $[a,b]$ and $\psi$ a smooth positive decreasing function defined on the same
interval. Suppose that $\phi'$ is monotonic in $[a,b]$ and that $|\phi'(s)|\ge\lambda$
for all $s\in[a,b]$. Then there is a universal constant $C>0$ such that
    \[
    \left|\int _a^b\e^{\mathrm{i}\phi(s)}\psi(s)\de s\right|\le C\frac{\psi(a)}{\lambda}.
    \]
\end{lemma}
A proof of this can be found in \cite{S}.
\end{section}

\begin{section} {Boundedness of $T_\nu^1$}
We need the following version of Schur's lemma.
\begin{lemma}
\label{schur}
Given two $\sigma$-finite measure spaces $(X,\,\mu)$, $(Y,\,\nu)$ and a
$\mu\otimes\nu$-measurable function $k$ on $X\times Y$, suppose that there exists
a positive constant $C$ such that
\[
\sup_{u\in X}\int_Y\int_X|k(x,\,y)k(u,\,y)|\de \mu(x)\de \nu(y)<C.
\]
Then, if $f\in L^2(X,\,\mu)$, the integral
\[
Kf(y)=\int_Xk(x,\,y)f(x)\de \mu(x)
\]
converges absolutely for $a.e.\,\, y\in Y $, the function $Kf$ thus defined is
in $L^2(Y,\,\nu)$ and
\[
\|Kf\|_2^2\le C\|f\|_2^2.
\]
\end{lemma}

The next proposition discusses the boundedness of the operator $T_\nu^1$.
\begin{proposition}
There is a positive constant $C$ such that for all
$\nu\ge 1$, for all intervals $I$, for all measurable functions $t(r)$ and
for all $g\in L^2(I)$,
\[
\|T_\nu^1 g\|_{L^2({0,\,1})}\le C\|g\|_{L^2(I)}.
\]
\end{proposition}
\noindent{\bf proof.}
The kernel of the operator $T_\nu^1$ is
\[
k(s,\,r)=\e^{\mathrm{i}t(r)s^2}\widetilde J_\nu(sr)\chi_{[0,\,\nu]}(sr)s^{-1/4},
\]
so that
$|k(s,\,r)|=|\widetilde J_\nu(sr)|\chi_{[0,\,\nu]}(sr)s^{-1/4}.$ By
Lemma \ref{schur},
\begin{eqnarray*}
\|T_\nu^1\|_2^2&\le&\sup_{u\in I}\int_0^1\widetilde J_\nu(uy)\chi_{[0,\,\nu]}(uy)u^{-1/4}
\int_I\widetilde J_\nu(sy)\chi_{[0,\,\nu]}(sy)s^{-1/4}\de s\de y\\
&=&\sup_{u\in I}\int_0^1\widetilde J_\nu(uy)\chi_{[0,\,\nu]}(uy)u^{-1/4}y^{-3/4}
\int_{I'\cap [0,\,\nu]}\widetilde J_\nu(t)t^{-1/4}\de t\de y.
\end{eqnarray*}
Since, by the well-known estimates for $J_\nu$ in the interval $[0,\,\nu/2]$
(see \cite{Wat}) and Stirling's formula,
\begin{eqnarray*}
\int_0^{\nu/2}\frac{J_\nu(t)}{t^\gamma}\de t &\le&
\int_0^{\nu/2}\frac{t^{\nu-\gamma}}{2^\nu\Gamma(\nu+1)}\de t=
\frac{\nu^{\nu+1-\gamma}}{2^{2\nu-\gamma+1}\Gamma(\nu+1)(\nu+1-\gamma)}\\
&\le&\frac C{\nu^{1/2+\gamma}}\left(\frac \e 4\right)^\nu,
\end{eqnarray*}
we have, using the estimate $\int_0^\nu|J_\nu(s)|\de s\le C$ (see \cite{GS}, Lemma 2.4),
\[
\int_0^\nu\frac{\widetilde J_\nu(t)}{t^{1/4}}\de t=
\int_0^{\nu/2}t^{1/4}{J_\nu(t)}\de t+
\int_{\nu/2}^\nu t^{1/4}{J_\nu(t)}\de t\le
C\nu^{1/4}.
\]
Therefore,
\begin{eqnarray*}
\|T_\nu^1\|_2^2&\le& C\sup_{u\in I}\frac{\nu^{1/4}}{u^{1/4}}
\int_0^1\frac{\widetilde J_\nu(uy)\chi_{[0,\,\nu]}(uy)}{y^{3/4}}\de y\\
&=&C\sup_{u\in I}\frac{\nu^{1/4}}{u^{1/2}}
\int_{[0,\,u]\cap[0,\,\nu]}J_\nu(s)s^{-1/4}\de s.
\end{eqnarray*}
Assume first that $u>\nu/2$. Then
\[
\frac{\nu^{1/4}}{u^{1/2}}
\int_{[0,\,u]\cap[0,\,\nu]}\frac{J_\nu(s)}{s^{1/4}}\de s\le
\frac C{\nu^{1/4}}\left[\int_0^{\nu/2}+\int_{\nu/2}^\nu\frac{ J_\nu(s)}
{s^{1/4}}\de s\right]\le \frac C{\nu^{1/2}}\le C.
\]
If instead $0<u<\nu/2$, then
\begin{eqnarray*}
&&\frac{\nu^{1/4}}{u^{1/2}}
\int_{[0,\,u]\cap[0,\,\nu]}\frac{J_\nu(s)}{s^{1/4}}\de s=
\frac{\nu^{1/4}}{u^{1/2}}
\int_0^u \frac{J_\nu(s)}{s^{1/4}}\de s\\
&\le&\frac{\nu^{1/4}}{u^{1/2}2^\nu\Gamma(\nu+1)}\int
_0^u s^{\nu-1/4}ds
\le C\frac{u^{\nu+1/4}\nu^{-3/4}\e^\nu}{2^\nu\nu^{\nu+1/2}}\\
&\le& C\frac{\nu^{\nu+1/4}\nu^{-3/4}\e^\nu}{4^\nu\nu^{\nu+1/2}}
=\frac C\nu\left(\frac e4\right)^\nu\le C.
\end{eqnarray*}
Thus  $\|T_\nu^1\|_2^2\le C$.
\qed

It is worth noting that in the study of $T_\nu^1$ we have not used the oscillation
given by $\e^{\mathrm{i}t(r)s^2}$.
\end{section}
\begin{section} {Boundedness of $T_\nu^2$}
Here we will use the following estimates on the Bessel functions: there
exists a positive constant $C$ such that if $s\in[\nu,\,\nu+\nu^{1/3}]$ then
$|\widetilde J_\nu(s)|\le C\nu^{1/6}$, and if $s\in[\nu+\nu^{1/3},\,2\nu]$ then
\[
|\widetilde J_\nu(s)|\le C\frac{\nu^{1/4}}{(s-\nu)^{1/4}}.
\]
These estimates are classical, but can be easily obtained from Lemma \ref{barcelo}
too. We can now state the boundedness result for $T_\nu^2$.

\begin{proposition}
There exists a positive constant $C$ such that
for all $\nu\ge 1$, for all intervals $I$, for all functions $t(r)$ and for all
$g\in L^2(I)$, we have
\[
\|T_\nu^2g\|_{L^2([0,\,1])}\le C\|g\|_{L^2(I)}.
\]
\end{proposition}
\noindent{\bf Proof.}
The absolute value of the kernel of $T_\nu^2$ is
\[
|k(s,\,r)|=|\widetilde J_\nu(sr)|\chi_{[\nu,\,\nu+\nu^{2/3}]}(sr)s^{-1/4}.
\]
By Schur's lemma,
\[
\|T_\nu^2\|_2^2\le\sup_{u\in I}\int_0^1|\widetilde J_\nu(uy)|\chi_{[\nu,\,\nu+\nu^{2/3}]}
(uy)u^{-1/4}\int_I|\widetilde J_\nu(sy)|\chi_{[\nu,\,\nu+\nu^{2/3}]}(sy)
s^{-1/4}\de s \de y.
\]
The innermost integral is bounded above by
\begin{eqnarray*}
&&\frac 1{y^{3/4}}\int_\nu^{\nu+\nu^{2/3}}|\widetilde J_\nu(t)|
t^{-1/4}dt\le\frac 1{y^{3/4}\nu^{1/4}}\int_\nu^{\nu+\nu^{2/3}}
|\widetilde J_\nu(t)|\de t\\
&\le&\frac C{y^{3/4}\nu^{1/4}}\left[\nu^{1/2}+
\int_{\nu+\nu^{1/3}}^{\nu+\nu^{2/3}}\frac{\nu^{1/4}}{(t-\nu)^{1/4}}\de t\right]\\
&\le& \frac C{y^{3/4}\nu^{1/4}}[\nu^{1/2}+\nu^{3/4}]
\le C\frac {\nu^{1/2}}{y^{3/4}}.
\end{eqnarray*}
Thus
\begin{eqnarray*}
\|T_\nu^2\|_2^2&\le& C\nu^{1/2}\sup_{u\in I}\int_0^1|\widetilde J_\nu(uy)|
\chi_{[\nu,\,\nu+\nu^{2/3}]}(uy)u^{-1/4}y^{-3/4}\de y\\
&=&C\nu^{1/2}\sup_{u\in I}u^{-1/2}\int_{[0,\,u]\cap[\nu,\,\nu+\nu^{2/3}]}
|\widetilde J_\nu(t)|t^{-3/4}\de t\\
&\le& C\nu^{-3/4}\int_\nu^{\nu+\nu^{2/3}}|\widetilde J_\nu(t)|\de t\\
&\le& C\nu^{-3/4+3/4}=C.
\end{eqnarray*}
\qed

Once more, in this proof we have not used the oscillation given by the exponential
nor the one given by the Bessel function.
\end{section}

\begin{section}{Boundedness of $T_\nu^3$.}
\begin{proposition}
There exists a positive constant $C$ such that
for all $\nu\ge 1$, for all intervals $I$, for all functions $t(r)$ and for all
$g\in L^2(I)$, we have
\[
\|T_\nu^3g\|_{L^2([0,\,1])}\le C\|g\|_{L^2(I)}.
\]
\end{proposition}
\noindent{\bf Proof.}
A trivial application of Cauchy-Schwartz's inequality yields
\begin{eqnarray*}
\| T_\nu^3g \|_{L^2([0,\,1])}^2 &=&
\int_0^1\left|\int_I\e^{\mathrm{i}t(r)s^2}r^{1/2}s^{1/4}
h_\nu(rs)\chi_{[\nu+\nu^{2/3},\,2\nu]}(rs)g(s)\de s\right|^2\de r\\
&\le& \int
_0^1\int_I|h_\nu(rs)|^2rs^{1/2}\chi_{[\nu+\nu^{2/3},\,2\nu]}(rs)
\de s\de r\,\|g\|_2^2\\
&\le&\int _0^1\int_{I'}|h_\nu(v)|^2\frac
{v^{1/2}}{r^{1/2}}\,\chi_{[\nu+\nu^{2/3},\,2\nu]}(v)\de v\de r\,\|g\|_2^2\\
&\le& C \|g\|_2^2\int_{\nu+\nu^{2/3}}^{2\nu}|h_\nu(v)|^2
v^{1/2}\de v\\
&=& C\nu^{3/2}
\|g\|_2^2\int_{1+\nu^{-1/3}}^{2}|h_\nu(\nu u)|^2
u^{1/2}\de u.
\end{eqnarray*}
The estimate
\[
|h_\nu(\nu u)|^2\le C\left(\frac{1}{\nu^3(u^2-1)^{7/2}} + \frac
2{\nu^{5/2}(u^2-1)^{7/4}u}+\frac 1{\nu^2u^2} \right),
\]
that holds for $u\in[1+\nu^{-2/3},\,2]$, concludes the proof.
\qed

\end{section}

\begin{section}{Boundedness of $T_\nu^4.$}
We shall need the following technical lemma. Its proof is a simple application
of the fundamental theorem of calculus.
\begin{lemma}
\label{technical}
Let $I$ be an interval and $g\in{\mathcal C}^3(I)$ be such that $g'(u)\le 0$,
$g''(u)\ge 0$ and $g'''(u)\le 0$ for all $u\in I$. Then for any $u,\,u_0\in I$,
\begin{enumerate}
\item if $u<u_0$, then $g(u)-g(u_0)\ge-g'(u_0)(u_0-u)$, and
\item if $u>u_0$, then $g(u_0)-g(u)\ge-g'(u_0)(u-u_0)-\frac 12 g''(u_0)(u-u_0)^2.$
\end{enumerate}
\end{lemma}
\begin{proposition}
There exists a positive constant $C$ such that
for all $\nu\ge 1$, for all intervals $I$, for all functions $t(r)$ and for all
$g\in L^2(I)$, we have
\[
\|T_\nu^4g\|_{L^2([0,\,1])}\le C\|g\|_{L^2(I)}.
\]
\end{proposition}
\noindent{\bf Proof.}
First write $T_\nu^4$ as the sum of two operators, by means of the equality
$\cos\theta=(\e^{\mathrm{i}\theta}+\e^{-\mathrm{i}\theta})/2$,
\begin{eqnarray*}
T_\nu^4g(r)&=&\sqrt{\frac 1{2\pi}}
\int_I \e^{\mathrm{i}t(r)s^2}\frac{r^{1/2}s^{1/4}\e^{\mathrm{i}\theta(rs)}}
{(s^2r^2-\nu^2)^{1/4}}
\chi_{[\nu+\nu^{2/3},\,2\nu]}(rs)
g(s)\de s+\\
&&+ \sqrt{\frac 1{2\pi}}
\int_I \e^{\mathrm{i}t(r)s^2}\frac{r^{1/2}s^{1/4}\e^{-\mathrm{i}\theta(rs)}}
{(s^2r^2-\nu^2)^{1/4}}
\chi_{[\nu+\nu^{2/3},\,2\nu]}(rs)
g(s)\de s.\\
\end{eqnarray*}
Observe that it is enough
to study just one of the two above operators, as long as we obtain
a result independent of the function $t(r)$, positive or negative. Let us then fix
our attention on the one with the $+$ sign
in the exponential (call it just $T$). The operator $TT^\ast$ has kernel
\[
K(r,\,\rho)=
\int_I\frac{\e^{\mathrm{i}[(t(r)-t(\rho))s^2+
\theta(rs)-\theta(\rho s)]}(r \rho s)^{1/2}\chi_{[\nu+\nu^{2/3},\,2\nu]}(rs)
\chi_{[\nu+\nu^{2/3},\,2\nu]}(\rho s)}
{(r^2s^2-\nu^2)^{1/4}(\rho^2s^2-\nu^2)^{1/4}}\de s.
\]
Let
\[
\tilde \theta (x)=\theta (\nu
x)=\nu\sqrt{x^2-1}-\nu\arccos(1/x)-\pi/4,\qquad x>1.
\]
Assuming $\rho<r$, calling $q=r/\rho$ and changing variables, $s=\nu
u/\rho$, we have that the kernel $K(r,\,\rho)$ equals
\[
\frac {\rho^{\beta-1/2}}{(r-\rho)^{\beta}}\left[\nu^{1/2}(q-1)^{\beta}
\int_{I\cap [1+\nu^{-1/3},\,2/q]}
\frac{\e^{\mathrm{i}[-a\nu u^2/2+\tilde\theta(qu)-\tilde\theta(u)]}}
{u^{1/2}(1-u^{-2})^{1/4}(1-q^{-2}u^{-2})^{1/4}}\de u\right],
\]
where $a=-2\nu(t(r)-t(\rho))/\rho^2$ and $\beta\in [1/2,\,1)$ will be fixed at our
convenience ($\beta=3/4$ will do).
Since the function $\min(r,\,\rho)^{\beta-1/2}|r-\rho|^{-\beta}$
is integrable in $\rho\in [0,\,1]$, uniformly in $r\in[0,\,1]$, by
Schur's lemma it is enough to show that the expression within
brackets is uniformly bounded in $a\in \Rset$, $\nu\gg 1$, $I$ any
interval, and $q\in(1,\,2)$ (for $q\ge 2$, the interval of integration becomes empty).

We introduce now some notation: for $u>1$ call
\begin{eqnarray*}
\psi(u)&=&\frac{\nu^{1/2}(q-1)^\beta}{u^{1/2}(1-u^{-2})^{1/4}(1-q^{-2}u^{-2})^{1/4}}
=\frac{\nu^{1/2}(q-1)^\beta u^{1/2}q^{1/2}}{(u^2-1)^{1/4}(q^2u^2-1)^{1/4}},\\
\phi(u)&=&-a\nu u^2/2+\tilde\theta(qu)-\tilde\theta(u),\\
\eta&=&-\log _\nu(q-1),
\end{eqnarray*}
so that $q=1+\nu^{-\eta}$, and the required uniformity in $q\in(1,\,2)$
is moved to the same one for $\eta> 0$.
Next observe that for $\eta\ge 1/(2\beta)$, the result is easily obtained since
\[
\left|\int_{I\cap[1+\nu^{-1/3},\,2/q]}\e^{\mathrm{i}\phi(u)}\psi(u)\de u\right|\le
\int_1^2\psi(u)\de u\le C\int_1^2\frac{\nu^{1/2-\eta\beta}}
{(u-1)^{1/2}}\de u\le C.
\]
Let us assume then that $0<\eta<1/(2\beta)$. This is the point where we start
using the oscillatory term in the estimation of our integral.
Since we want to use Van der Corput's lemma, we need to study the
function $\phi'$. Note that
\[
\phi'(u)=\nu\left(\sqrt{q^2-u^{-2}}-\sqrt{1-u^{-2}}-au\right)=\nu(f(u)-au),
\]
where $f$ is implicitly  defined by the above equality.
Let us begin by considering only those values of $a$ for which there is
a zero of $\phi'$ in the interval
$[1+\nu^{-1/3},\,2/q]$. Thus, parametrize
$a$ in such a way that this zero is $1+\nu^{-\gamma}$, with
$\gamma\in[0,\,1/3]$. This way, the required uniformity in the parameter $a$ is moved
to the parameter $\gamma$. For further reference, observe that
\[
a=\frac{\sqrt{q^2(1+\nu^{-\gamma})^2-1}-\sqrt{(1+\nu^{-\gamma})^2-1}}
{(1+\nu^{-\gamma})^2}.
\]
\begin{figure}[htb!]
\begin{center}
  \includegraphics[width=12cm,height=9cm]{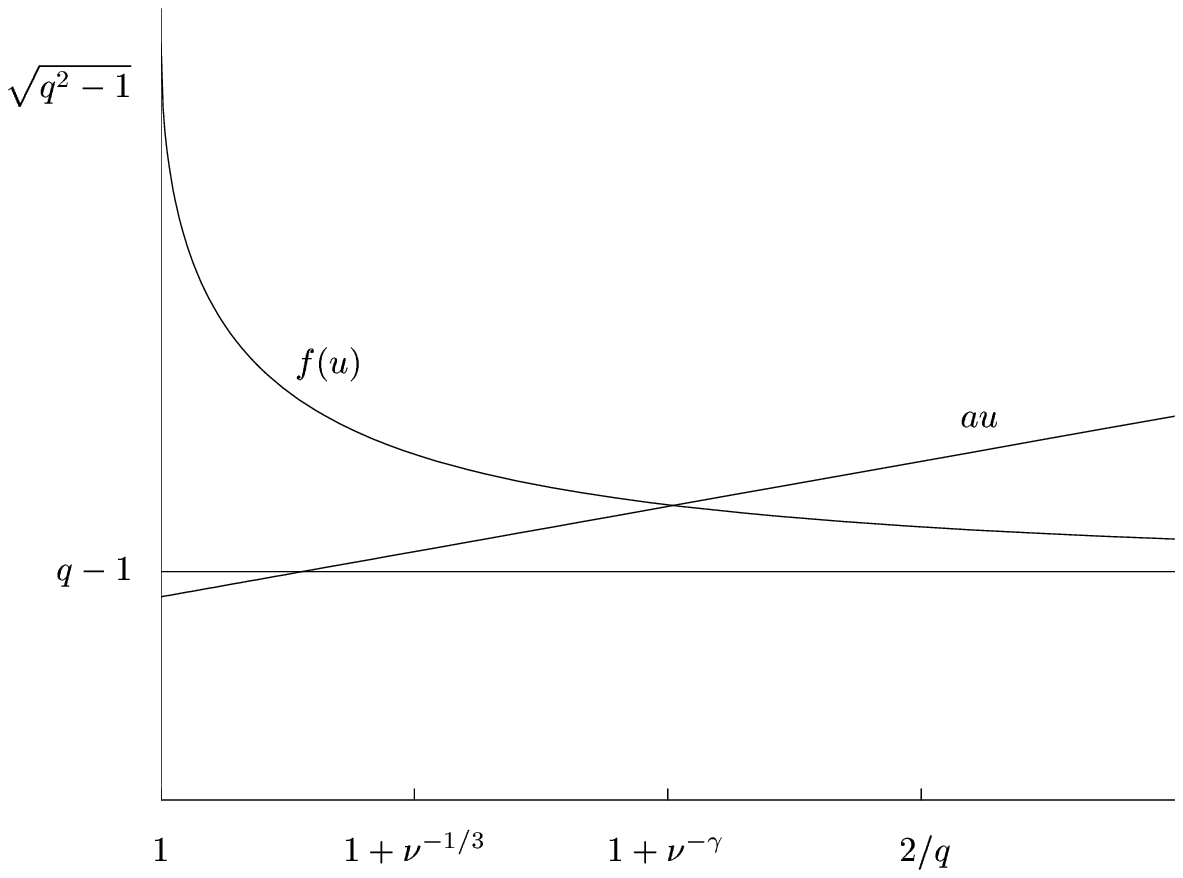}
  \caption{\label{Fm}The curves $f(u)$ and $au$.}
\end{center}
\end{figure}
Let $u\in[1+\nu^{-1/3},\,2/q]$.
One can easily see that $f$ satisfies all the hypotheses of
Lemma~\ref{technical}. Thus, recalling that $\phi'(u)=
\nu(f(u)-au)$, we may deduce that if $u < 1+\nu^{-\gamma}$ then
\begin{equation}
\label{derivativeone}
|\phi'(u)|\ge \nu(1+\nu^{-\gamma}-u)\left(a-f'(1+\nu^{-\gamma})\right),
\end{equation}
whereas if $u>1+\nu^{-\gamma}$
\begin{equation}
\label{derivativetwo}
|\phi'(u)|\ge \nu(u-1-\nu^{-\gamma})\left(a-f'(1+\nu^{-\gamma})-
\frac 12 f''(1+\nu^{-\gamma})(u-1-\nu^{-\gamma})\right).
\end{equation}

Define $\xi=\min(\eta,\,\gamma)$ and
\[
\delta=\frac 12-\eta\beta+\frac\gamma 4+\frac\xi 4.
\]
Observe that $(\eta,\,\gamma)$ may vary in the rectangle ${\mathcal R}=(0,\,1/(2\beta))
\times[0,\,1/3]$. Divide $\mathcal R$ into two regions, ${\mathcal F}=\{(\eta,\,\gamma)
\in\mathcal R
\colon \delta\ge\gamma\}$ and $\mathcal G=\mathcal R\setminus \mathcal F$.
\begin{figure}[htb!]
\begin{center}
  \includegraphics[width=12cm,height=9cm]{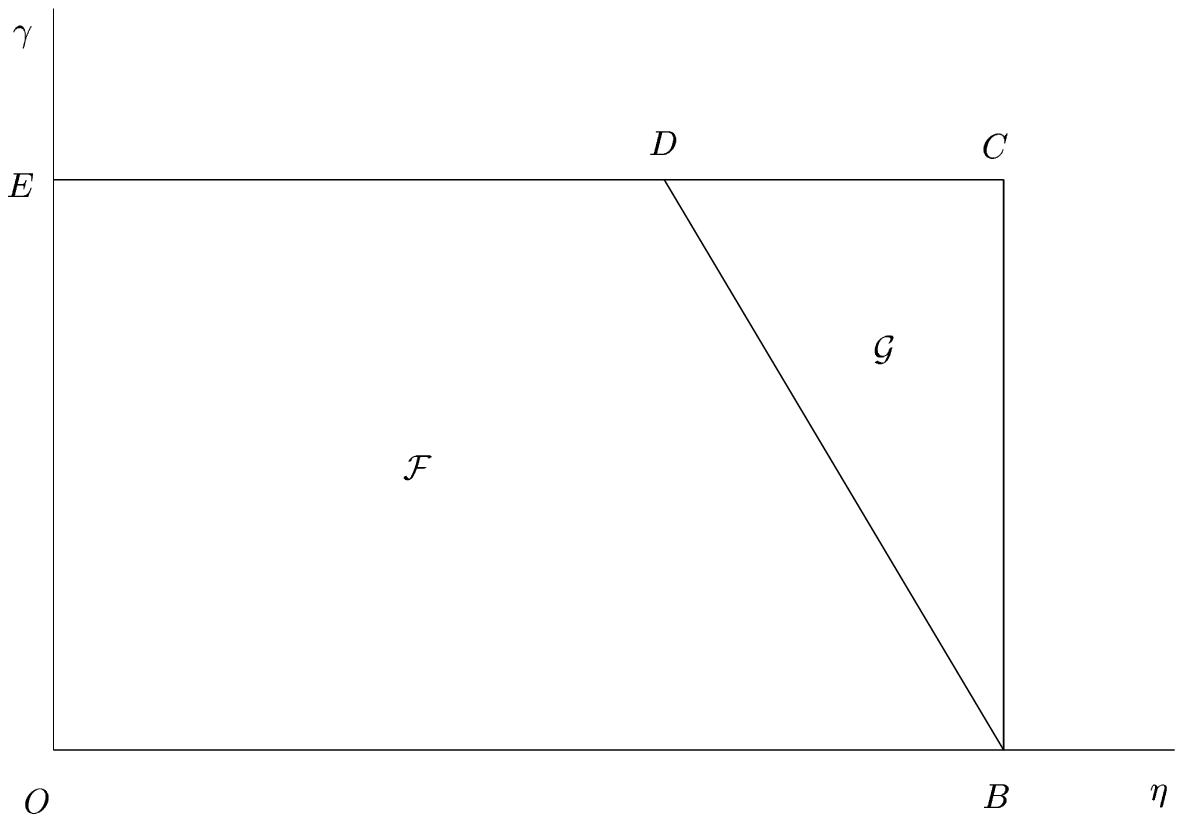}
  \caption{\label{FL} The sets $\mathcal F$ and $\mathcal G$, where $B=(1/(2\beta),\,0)$,
  $C=(1/(2\beta),\,1/3)$, $D=(1/(3\beta),\,1/3)$ and $E=(0,\,1/3)$.}
\end{center}
\end{figure}

Consider first the case $(\eta,\,\gamma)\in\mathcal F$.
Divide the interval $[1+\nu^{-1/3},\,2/q]$ into the union of four subintervals
(defined to be empty when the left endpoint happens to be bigger than the
right endpoint):
\begin{eqnarray*}
{\mathcal A}_1&=&[1+\nu^{-1/3},\,1+\nu^{-\gamma}/10],\\
{\mathcal A}_2&=&[1+\nu^{-\gamma}/10,\,1+\nu^{-\gamma}-\nu^{-\delta}/N],\\
{\mathcal A}_3&=&[1+\nu^{-\gamma}-\nu^{-\delta}/N,\,
          1+\nu^{-\gamma}+\nu^{-\delta}/N],\\
{\mathcal A}_4&=&[1+\nu^{-\gamma}+\nu^{-\delta}/N,\,2/q],\\
\end{eqnarray*}
where $N$ is a large number that will be fixed at our convenience.
The interval $\mathcal A_3$ is a neighborhood of the zero of $\phi'$, where the
oscillation vanishes. The best we can do here is then to estimate the corresponding
integral with the magnitude of the integrand:
\begin{eqnarray*}
\left|\int_{I\cap{\mathcal A}_3}\e^{\mathrm{i}\phi(u)}\psi(u)\de u\right|
&\le&\int_{1+\nu^{-\gamma}-\nu^{-\delta}/N}^{1+\nu^{-\gamma}+\nu^{-\delta}/N}
\psi(u)\de u\le \frac 2{N\nu^{\delta}}\,\psi\left(1+\frac 1{2\nu^{\gamma}}\right)\\
&\le& C\frac 1{\nu^{\delta}}\frac{\nu^{1/2-\eta\beta}}{\nu^{-\gamma/4}((1+\nu^{-\eta})
(1+\nu^{-\gamma}/2)-1)^{1/4}}\\
&\le&
C\frac 1{\nu^{\delta}}\frac{\nu^{1/2-\eta\beta+\gamma/4}}{(\nu^{-\eta}
+\nu^{-\gamma})^{1/4}}\le C\frac{\nu^\delta}{\nu^\delta}\le C.\\
\end{eqnarray*}
For $\mathcal A_1$, we can use Van der Corput's lemma. Thus
\[
\left|\int_{I\cap{\mathcal A_1}}\e^{\mathrm{i}\phi(u)}\psi(u)\de u\right|
\le C\,\frac{\psi(1+\nu^{-1/3})}{\phi'(1+\nu^{-\gamma}/10)}.
\]
Observe that
\[
\psi(1+\nu^{-1/3})\le C\frac {\nu^{1/2-\eta\beta}}
{\nu^{-1/12}(\nu^{-1/3}+\nu^{-\eta})^{1/4}}\le C\,\nu^{7/12-\eta\beta+\zeta/4},
\]
where $\zeta=\min(\eta,\,1/3)$. As for $\phi'$, using
(\ref{derivativeone}) we have that
\begin{eqnarray*}
|\phi'(1+\nu^{-\gamma}/10)|&\ge& \nu(1+\nu^{-\gamma}-1-\nu^{-\gamma}/10)
\left(a-f'(1+\nu^{-\gamma})\right)\\
&\ge&\frac 9{10}\,\nu^{1-\gamma}
\left(-f'(1+\nu^{-\gamma})\right)\\
&=&\frac 9{10}\,\nu^{1-\gamma}\left(
-\frac{(1+\nu^{-\gamma})^{-2}}{\sqrt{q^2(1+\nu^{-\gamma})^2-1}}
+\frac{(1+\nu^{-\gamma})^{-2}}{\sqrt{(1+\nu^{-\gamma})^2-1}}
\right)\\
&\ge&C\,\nu^{1-\gamma}\frac{\sqrt{(1+\nu^{-\eta})^2(1+\nu^{-\gamma})^2-1}
-\sqrt{(1+\nu^{-\gamma})^2-1}}{\sqrt{(1+\nu^{-\eta})^2(1+\nu^{-\gamma})^2-1}
\sqrt{(1+\nu^{-\gamma})^2-1}}\\
&\ge& C\nu^{1-\gamma/2-\eta+\xi}.\\
\end{eqnarray*}
Thus
\begin{eqnarray*}
\left|\int_{I\cap{\mathcal A_1}}\e^{\mathrm{i}\phi(u)}\psi(u)\de u\right|
&\le& C\frac{\nu^{7/12-\eta\beta+\zeta/4}}{\nu^{1-\gamma/2-\eta+\xi}}\\
&\le& C\nu^{-5/12+\eta(1-\beta)+\gamma/2+\zeta/4-\xi}
\le C\nu^{-1/6+(1-\beta)/(2\beta)}\le C,
\end{eqnarray*}
if $\beta\ge 3/4$.

Let us now consider $\mathcal A_2$. Once again, using Van der Corput's lemma,
\[
\left|\int_{{\mathcal A}_2\cap I}\e^{\mathrm{i}\phi(u)}\psi(u)\de u\right|
\le C\frac{\psi(1+\nu^{-\gamma}/10)}{\phi'(1+\nu^{-\gamma}-\nu^{-\delta}/N)}.
\]
Proceeding as in the previous case, we see that
\[
\psi(1+\nu^{-\gamma}/10)\le C\frac {\nu^{1/2-\eta\beta}}
{\nu^{-\gamma/4}(\nu^{-\gamma}+\nu^{-\eta})^{1/4}}
\le C\,\nu^{1/2-\eta\beta+\gamma/4+\xi/4},
\]
and that
\begin{eqnarray*}
|\phi'(1+\nu^{-\gamma}-\nu^{-\delta}/N)|&\ge& \nu(1+\nu^{-\gamma}-1-
\nu^{-\gamma}+\nu^{-\delta}/N)
\left(a-f'(1+\nu^{-\gamma})\right)\\
&\ge&\frac{\nu^{1-\delta}}{N}
\left(-f'(1+\nu^{-\gamma})\right)\\
&\ge&C\nu^{1-\delta-\eta+\gamma/2+\xi}\\
&=& C\nu^{1/2-\eta(1-\beta)+\gamma/4+3\xi/4}.\\
\end{eqnarray*}
Therefore
\begin{eqnarray*}
\left|\int_{I\cap{\mathcal A_2}}\e^{\mathrm{i}\phi(u)}\psi(u)\de u\right|
&\le& C\frac{\nu^{1/2-\eta\beta+\gamma/4+\xi/4}}
{\nu^{1/2-\eta(1-\beta)+\gamma/4+3\xi/4}}\\
&\le& C\nu^{-(2\beta-1)\eta-\xi/2}\le C,
\end{eqnarray*}
if $\beta\ge 1/2.$

Let us now move to the study of the interval $\mathcal A_4$.
Using Van der Corput's lemma, we have
\[
\left|\int_{{\mathcal A}_4\cap I}\e^{\mathrm{i}\phi(u)}\psi(u)\de u\right|
\le C\frac{\psi(1+\nu^{-\gamma}+\nu^{-\delta}/N)}
{|\phi'(1+\nu^{-\gamma}+\nu^{-\delta}/N)|}.
\]
As usual, we see that
\[
\psi(1+\nu^{-\gamma}+\nu^{-\delta}/N)\le C\frac {\nu^{1/2-\eta\beta}}
{\nu^{-\gamma/4}(\nu^{-\gamma}+\nu^{-\eta})^{1/4}}
\le C\,\nu^{1/2-\eta\beta+\gamma/4+\xi/4},
\]
while using (\ref{derivativetwo}), we see that
\begin{eqnarray*}
|\phi'(1+\nu^{-\gamma}+\nu^{-\delta}/N)|&\ge& \nu\,\frac 1{N\nu^{\delta}}
\left(a-f'(1+\nu^{-\gamma})-\frac 12f''(1+\nu^{-\gamma})\frac 1 {N\nu^\delta}\right)\\
&\ge&\frac{\nu^{1-\delta}}{N}
\left(-f'(1+\nu^{-\gamma})-\frac 12f''(1+\nu^{-\gamma})
\frac 1{N\nu^{\delta}}\right).\\
\end{eqnarray*}
We already have the estimate $-f'(1+\nu^{-\gamma})\ge C\nu^{\gamma/2+\xi-\eta}.$
We shall now show that there is a positive constant $C$ such that
\begin{equation}
\label{secondderivative}
|f''(1+\nu^{-\gamma})|\le C\nu^{3\gamma/2+\xi-\eta}.
\end{equation}
In the following computations, we will call $u_0=1+x=1+\nu^{-\gamma}$,
 $q=1+y=1+\nu^{-\eta}$, with $x,\,y\in[0,\,1]$, and $z=\max(x,\,y)$. Thus
\begin{eqnarray*}
&&|f''(u_0)|=\frac{(3u_0^2-2)(q^2u_0^2-1)^{3/2}-(3q^2u_0^2-2)(u_0^2-1)^{3/2}}
{u_0^3(u_0^2-1)^{3/2}(q^2u_0^2-1)^{3/2}}=\\
&&\frac{(3u_0^2-2)^2(q^2u_0^2-1)^{3}-(3q^2u_0^2-2)^2(u_0^2-1)^{3}}
{u_0^3(u_0^2-1)^{3/2}(q^2u_0^2-1)^{3/2}\left(
(3u_0^2-2)(q^2u_0^2-1)^{3/2}+(3q^2u_0^2-2)(u_0^2-1)^{3/2}\right)}.\\
\end{eqnarray*}
The numerator of the above expression is a polynomial in $x$ and $y$, sum of monomials
of degrees $3$ to $16$, none of which is of the form $x^j$ for any $j$.
Therefore this numerator is bounded above in absolute value by
\[
Cy(x^2+xy+y^2)\le Cyz^2.
\]
On the other hand, the denominator is bounded below
in absolute value by
\[
Cx^{3/2}(x+y)^{3/2}((x+y)^{3/2}+x^{3/2})\ge Cx^{3/2}z^3.
\]
It follows that
\[
|f''(u_0)|\le C\frac{yz^2}{x^{3/2}z^3}=C\nu^{3\gamma/2+\xi-\eta},
\]
as desired.
Thus we may deduce that
\begin{eqnarray*}
|\phi'(1+\nu^{-\gamma}+\nu^{-\delta}/N)|&\ge& \frac CN\nu^{1-\delta}
\left(\nu^{\gamma/2+\xi-\eta}- \frac {C\nu^{3\gamma/2+\xi-\eta}}{N\nu^\delta}\right)\\
&\ge&\frac CN\nu^{1-\delta}\nu^{\gamma/2+\xi-\eta}(1-\frac CN\nu^{\gamma-\delta} )\\
&\ge& C\nu^{1-\delta+\gamma/2+\xi-\eta}\\
&\ge& C\nu^{1/2+\gamma/4+3\xi/4-\eta(1-\beta)},\\
\end{eqnarray*}
if we take $N$ big enough (recall we are in the case $\delta\ge\gamma$). We may now
conclude
\[
\left|\int_{{\mathcal A}_4\cap I}\e^{\mathrm{i}\phi(u)}\psi(u)\de u\right|
\le
C\frac{\nu^{1/2-\eta\beta+\gamma/4+\xi/4}}{\nu^{1/2+\gamma/4+3\xi/4-\eta(1-\beta)}}
\le C\nu^{-\eta(2\beta-1)-\xi/2}\le C,
\]
if $\beta\ge 1/2$.

It remains to study the case
$(\eta,\,\gamma)\in\mathcal G$, that is
$\delta<\gamma$. Observe that this implies $\gamma\le\eta$ and therefore $\xi=\gamma.$
Divide the interval $[1+\nu^{-1/3},\,2/q]$ into the union of three subintervals
(defined to be empty when the left endpoint happens to be bigger than the
right endpoint):
\begin{eqnarray*}
{\mathcal A}_1&=&[1+\nu^{-1/3},\,1+\nu^{-\gamma}/10],\\
{\mathcal A}_2&=&[1+\nu^{-\gamma}/10,\,1+2\nu^{-1/2+\eta\beta-\gamma/2}],\\
{\mathcal A}_3&=&[1+2\nu^{-1/2+\eta\beta-\gamma/2},\,2/q],\\
\end{eqnarray*}
The interval $\mathcal A_2$ is a neighborhood of the zero of $\phi'$, where the
oscillation vanishes, so we estimate the associated integral with the magnitude
of the integrand:
\begin{eqnarray*}
\left|\int_{I\cap{\mathcal A}_2}\e^{\mathrm{i}\phi(u)}\psi(u)\de u\right|
&\le&\int_{1+\nu^{-\gamma}/10}^{1+2\nu^{-1/2+\eta\beta-\gamma/2}}
\psi(u)\de u\\
&\le& 2\nu^{-1/2+\eta\beta-\gamma/2}\,
\psi\left(1+\nu^{-\gamma}/10\right)\\
&\le& C\frac{\nu^{-1/2+\eta\beta-\gamma/2}\nu^{1/2-\eta\beta}}
{\nu^{-\gamma/4}((1+\nu^{-\eta})
(1+\nu^{-\gamma}/10)-1)^{1/4}}\\
&\le&
C\frac{\nu^{-\gamma/2}}{\nu^{-\gamma/2}}\le C.\\
\end{eqnarray*}

The study of $\mathcal A_1$ is exactly the same as in the case $\gamma\le\delta$,
thus we do not repeat it.

As for $\mathcal A_3$, we use Van der Corput's lemma, obtaining
\[
\left|\int_{{\mathcal A}_3\cap I}\e^{\mathrm{i}\phi(u)}\psi(u)\de u\right|\le
C\frac{\psi(1+2\nu^{-1/2+\eta\beta-\gamma/2})}{|\phi'(1+2\nu^{-1/2+\eta\beta-\gamma/2})|}.
\]
Since $a$ is positive and $f$ is decreasing and $f(1+\nu^{-\gamma})=a(1+\nu^{-\gamma})$,
 we may say that
\begin{eqnarray*}
|\phi'(1+2\nu^{-1/2+\eta\beta-\gamma/2})|&\ge&
\nu(a(1+2\nu^{-1/2+\eta\beta-\gamma/2})-f(1+\nu^{-\gamma}))\\
&=&\nu a(2\nu^{-1/2+\eta\beta-\gamma/2}-\nu^{-\gamma})
\ge\nu a\nu^{-1/2+\eta\beta-\gamma/2}\\
&\ge&\nu^{1/2+\eta\beta-\gamma/2}\frac{q^2-1}{\sqrt{q^2(1+\nu^{-\gamma})^2-1}+
\sqrt{(1+\nu^{-\gamma})^2-1}}\\
&\ge& C\nu^{1/2+\eta\beta-\gamma/2}\frac{\nu^{-\eta}}{(\nu^{-\gamma}+\nu^{-\eta})^{1/2}
+\nu^{-\gamma/2}}\\
&\ge& C\nu^{1/2+\eta\beta-\gamma/2}\nu^{-\eta+\gamma/2}=C\nu^{1/2+\eta\beta-\eta}.\\
\end{eqnarray*}
On the other hand,
\begin{eqnarray*}
\psi(1+2\nu^{-1/2+\eta\beta-\gamma/2})&\le& C\frac{\nu^{1/2-\eta\beta}}
{\nu^{-1/8+\eta\beta/4-\gamma/8}(\nu^{-\eta}+\nu^{-1/2+\eta\beta-\gamma/2})^{1/4}}\\
&\le&
C\frac{\nu^{1/2-\eta\beta}}{\nu^{-1/4+\eta\beta/2-\gamma/4}}
=C\nu^{3/4-3\eta\beta/2+\gamma/4}.\\
\end{eqnarray*}
Therefore,
\[
\left|\int_{{\mathcal A}_3\cap I}\e^{\mathrm{i}\phi(u)}\psi(u)\de u\right|\le
C\,\frac{\nu^{3/4-3\eta\beta/2+\gamma/4}}{\nu^{1/2+\eta\beta-\eta}}=
C\nu^{1/4-5\eta\beta/2+\gamma/4+\eta}\le C,
\]
if $(\eta,\,\gamma)\in\mathcal G$, and $\beta\ge 2/3$.

It remains to study the boundedness of the integral
\[
\left|\int_{I\cap[1+\nu^{-1/3},\,2/q]}\e^{\mathrm{i}\phi(u)}\psi(u)\de u\right|
\]
for the values of $a$ for which $\phi'$ has no zeros in $[1+\nu^{-1/3},\,2/q]$.
Call $a_0$ and $a_1$ the values of $a$ for which the zero of $\phi'$ is
$1+\nu^{-1/3}$ and $2/q$, respectively. Geometrically, it is clear that
for any fixed $u\in[1+\nu^{-1/3},\,2/q]$, the value of $|\phi'(u)|$ grows as
$a$ goes from $a_0$ to $\infty$, and decreases as $a$ goes from $-\infty$ to $a_1$,
while $\psi(u)$ stays unchanged. Thus, all the estimates we obtained for
$a_0$ using Van der Corput's lemma or simply the magnitude of $\psi$, remain true
for any $a\ge a_0$, and those we obtained for $a_1$ remain true for any $a\le a_1$.
This concludes the proof.
\qed

\end{section}
\begin{section}{Boundedness of $T_\nu^5.$}
\begin{proposition}
There exists a positive constant $C$ such that
for all $\nu\ge 1$, for all intervals $I$, for all functions $t(r)$ and for all
$g\in L^2(I)$, we have
\[
\|T_\nu^5g\|_{L^2([0,\,1])}\le C\|g\|_{L^2(I)}.
\]
\end{proposition}
\noindent{\bf Proof.}
The kernel of the operator $T_\nu^5(T_\nu^5)^\ast$ is
\[
L(r,\,\rho)=\int_I\e^{\mathrm{i}(t(r)-t(\rho))s^2}\widetilde h_\nu(rs)
\widetilde h_\nu(\rho s)\chi_{[2\nu,\,\infty)}(rs)\chi_{[2\nu,\,\infty)}(\rho s)
s^{-1/2}\de s,
\]
Thus, using Lemma \ref{barcelo},
\[
|L(r,\,\rho)|\le\int_{\frac{2\nu}{\min(r,\,\rho)}}^{\infty}|\widetilde h_\nu(rs)
\widetilde h_\nu(\rho s)|s^{-1/2}\de s\le  C
\int_{\frac{2\nu}{\min(r,\,\rho)}}^{\infty}\frac {s^{-3/2}}{\sqrt{r\rho}}\de s\le
\frac{C\sqrt{\min(r,\,\rho)}}{\sqrt{\nu r\rho}}.
\]
Since
\[
\int_0^1|L(r,\,\rho)|\de r\le \int_0^1\frac{C\sqrt{\min(r,\,\rho)}}{\sqrt{\nu r\rho}}
\de r=\frac C{\sqrt{\nu\rho}}\int_0^\rho \de r+
\frac C{\sqrt{\nu}}\int_\rho^1\frac {\de r}
{\sqrt r}=\frac C{\sqrt \nu}(2-\sqrt\rho)
\]
is uniformly bounded in $\rho\in [0,\,1]$, by Schur's lemma the operators
$T_\nu^5(T_\nu^5)^\ast$ are uniformly bounded, and so are the $T_\nu^5$'s.
\qed

\end{section}

\begin{section}{Boundedness of $T_\nu^6.$}
\begin{proposition}
There exists a positive constant $C$ such that
for all $\nu\ge 1$, for all intervals $I$, for all functions $t(r)$ and for all
$g\in L^2(I)$, we have
\[
\|T_\nu^6g\|_{L^2([0,\,1])}\le C\|g\|_{L^2(I)}.
\]
\end{proposition}
Proceeding as for $T_\nu^4$, write $T_\nu^6$ as the sum of two operators,
by means of the equality
$\cos\theta=(\e^{\mathrm{i}\theta}+\e^{-\mathrm{i}\theta})/2$,
\begin{eqnarray*}
T_\nu^6g(r)&=&\sqrt{\frac 1{2\pi}}
\int_I \e^{\mathrm{i}t(r)s^2}\frac{r^{1/2}s^{1/4}\e^{\mathrm{i}\theta(rs)}}
{(s^2r^2-\nu^2)^{1/4}}
\chi_{[2\nu,\,\infty)}(rs)g(s)\de s+\\
&&+ \sqrt{\frac 1{2\pi}}
\int_I \e^{\mathrm{i}t(r)s^2}\frac{r^{1/2}s^{1/4}\e^{-\mathrm{i}\theta(rs)}}
{(s^2r^2-\nu^2)^{1/4}}
\chi_{[2\nu,\,\infty)}(rs)
g(s)\de s.\\
\end{eqnarray*}
Once again, it is enough
to study just one of these two operators, for example the one with the $+$ sign
in the exponential (call it just $T$). The operator $TT^\ast$ has kernel
\[
K(r,\,\rho)=
\int_I\frac{\e^{\mathrm{i}[(t(r)-t(\rho))s^2+
\theta(rs)-\theta(\rho s)]}r^{1/2}\rho^{1/2}s^{1/2}\chi_{[2\nu,\,\infty)}(rs)
\chi_{[2\nu,\,\infty)}(\rho s)}
{(r^2s^2-\nu^2)^{1/4}(\rho^2s^2-\nu^2)^{1/4}}\de s.
\]
Assuming $\rho<r$, calling $p=(r-\rho)/\rho$ and $\sigma=p\nu$,
and changing variables, $s=u/(r-\rho)$, we have the kernel
\[
\frac 1{(r-\rho)^{1/2}}
\int_{I\cap [2\sigma,\,\infty)}
\frac{\e^{\mathrm{i}[-a u^2/2+\theta((p+1)u/p)-\theta(u/p)]}}
{u^{1/2}(1-\sigma^2(p+1)^{-2}u^{-2})^{1/4}(1-\sigma^2u^{-2})^{1/4}}\de u,
\]
where $a=-2(t(r)-t(\rho))/(r-\rho)^2$. Since the function $|r-\rho|^{-1/2}$
is integrable in $r$, uniformly in $\rho$, by Schur's lemma it is enough
to show that the integral is uniformly bounded in the interval $I$,
in $p>0$, in $\sigma>0$, and in $a\in\Rset$.
Let us call
\begin{eqnarray*}
\phi(u)&=&-\frac a2 u^2+\theta\left(\frac{p+1}pu\right)-\theta\left(\frac up\right)\\
\psi(u)&=&\frac 1{u^{1/2}\left(1-\frac{\sigma^2}{(p+1)^{2}u^{2}}\right)^{1/4}
\left(1-\frac{\sigma^2}{u^2}\right)^{1/4}}.\\
\end{eqnarray*}
Observe that
\[
\phi'(u)=-au+\frac{(p+2)u}{\sqrt{(p+1)^2u^2-\sigma^2}+\sqrt{u^2-\sigma^2}}=-au+f(u),
\]
(note that here ``$f$'' indicates a different function from the one in section~5)
and that the function $\psi$ is decreasing with
\[
\psi(u)\le\frac 1{\sqrt{u-\sigma}}.
\]
Note that, since $\phi'$ is the difference of a concave up function and a linear
function, $\phi''$ is the difference of an increasing function and a constant.
Hence, $\phi''$ is increasing and therefore it changes sign at most once. By assuming
that the interval $I$ is contained in an interval where $\phi''$ has constant sign,
we can apply Van der Corput's lemma to
\[
\left|\int_{I\cap[2\sigma,\,\infty)}\e^{\mathrm{i}\phi(u)}\psi(u)\de u\right|.
\]
In order to do it, we need to study the function $\phi'$.
As usual, we consider only those values of $a$ for which there is
a zero of $\phi'$ in the interval
$[2\sigma,\,\infty)$, that is
\[
0<a\le\frac{\sqrt{4(p+1)^2-1}-\sqrt{3}}
{4\sigma p}.
\]
Assume first that $\sigma\ge 1$.
Let us parametrize
$a$ in such a way that the zero of $\phi'$ is $\sigma+\sigma^{\gamma}$, with
$\gamma\ge 1$. This gives
\[
a=\frac{(p+2)}{\sqrt{(p+1)^2(\sigma+\sigma^\gamma)^2-\sigma^2}
+\sqrt{(\sigma+\sigma^\gamma)^2-\sigma^2}}.
\]
In this way, the required uniformity in the parameter $a$ is equivalent
to the uniformity in the parameter $\gamma$.
\begin{figure}[htb!]
\begin{center}
  \includegraphics[width=12cm,height=9cm]{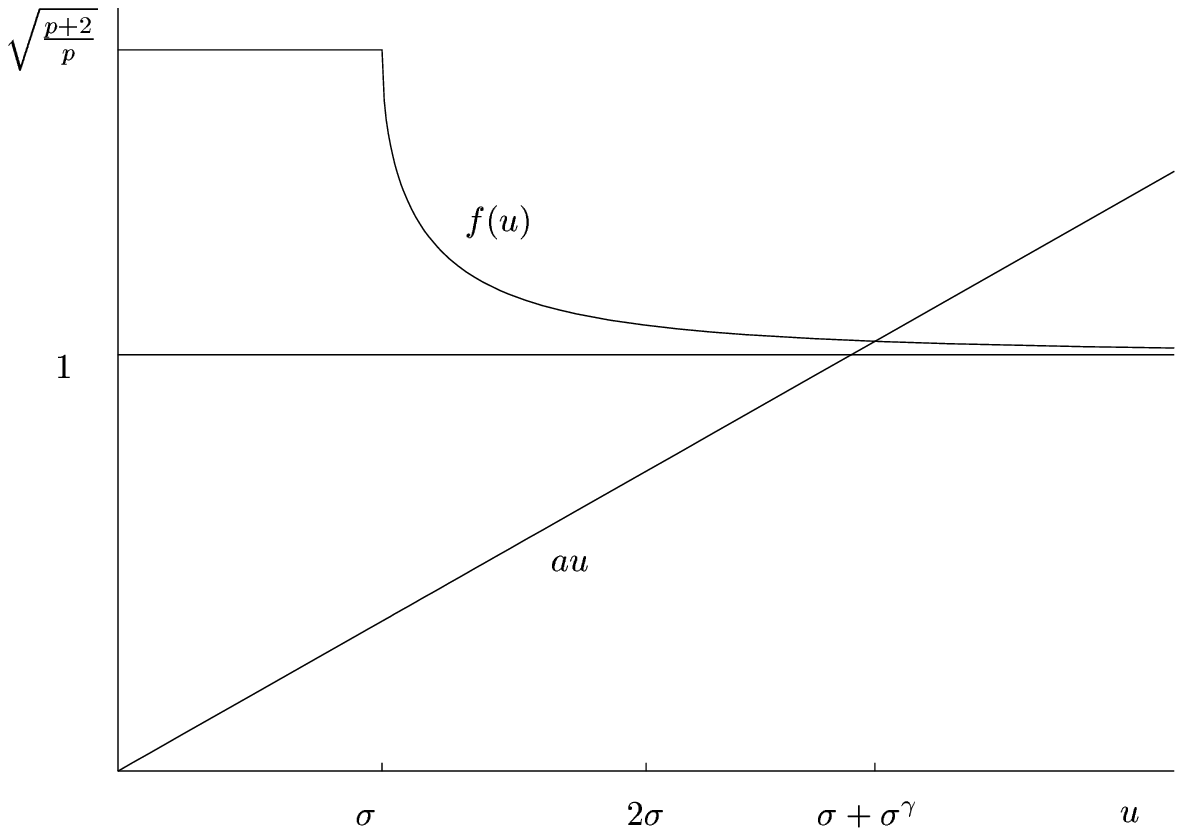}
  \caption{\label{FN}The curves $f(u)$ and $au$.}
\end{center}
\end{figure}
In order to apply Van der Corput's lemma we need to estimate $|\phi'|$ from below.
Observe that
\begin{eqnarray*}
|\phi'(u)|
&=&\left|-\frac{(p+2)u}{\sqrt{(p+1)^2(\sigma+\sigma^\gamma)^2-\sigma^2}
+\sqrt{(\sigma+\sigma^\gamma)^2-\sigma^2}}+\right.\\
&&+\left.\frac{(p+2)u}{\sqrt{(p+1)^2u^2-\sigma^2}+\sqrt{u^2-\sigma^2}}\right|\\
&=&\left|\frac{(p+2)u}{(\sqrt{(p+1)^2u^2-\sigma^2}+\sqrt{u^2-\sigma^2})}\times\right.\\
&&\times\frac 1
{({\sqrt{(p+1)^2(\sigma+\sigma^\gamma)^2-\sigma^2}
+\sqrt{(\sigma+\sigma^\gamma)^2-\sigma^2}})}\times\\
&&\times\left(
\frac{(p+1)^2((\sigma+\sigma^\gamma)^2-u^2)}
{\sqrt{(p+1)^2(\sigma+\sigma^\gamma)^2-\sigma^2}
+\sqrt{(p+1)^2u^2-\sigma^2}}+\right.\\
&&\left.\left.+\frac{(\sigma+\sigma^\gamma)^2-u^2}
{\sqrt{(\sigma+\sigma^\gamma)^2-\sigma^2}
+\sqrt{u^2-\sigma^2}}\right)\right|\\
&\ge&\left|\frac{(p+2)u}{((p+1)u+u)((p+1)(\sigma+\sigma^\gamma)+(\sigma+\sigma^\gamma))}
\times\right.\\
&&\times\left.\left(\frac{(p+1)^2((\sigma+\sigma^\gamma)^2-u^2)}
{(p+1)(\sigma+\sigma^\gamma)+(p+1)u}+
\frac{(\sigma+\sigma^\gamma)^2-u^2}{\sigma+\sigma^\gamma+u}\right)\right|\\
&=&\left|\frac{\sigma+\sigma^\gamma-u}{\sigma+\sigma^\gamma}\right|.
\end{eqnarray*}
Next divide the interval $[2\sigma,\,\infty)$ into four subintervals, given by the
following partition
\begin{eqnarray*}
u_1&=&2\sigma,\\
u_2&=&\max(2\sigma,\,\sigma+\sigma^\gamma/2),\\
u_3&=&\max(u_2,\,\sigma+\sigma^\gamma-\sigma^{\gamma/2}),\\
u_4&=&\sigma+\sigma^\gamma+\sigma^{\gamma/2},\\
\end{eqnarray*}
and study each case separately.
Applying Van der Corput's lemma and using the above estimates for $\psi$ and
$\phi'$, we obtain that when $[u_1,\,u_2]$ is non-degenerate,
\[
\left|\int_{[u_1,\,u_2]\cap I}\e^{\mathrm{i}\phi(u)}\psi(u)\de u\right|\le
C\frac{\psi(2\sigma)}{|\phi'(\sigma+\sigma^\gamma/2)|}
\le C\frac{\sigma+\sigma^\gamma}{\sqrt{\sigma}\sigma^{\gamma}/2}
\le \frac C{\sqrt\sigma}\le C.
\]
On the other hand, when $[u_2,\,u_3]$ is non-degenerate,
\[
\left|\int_{[u_2,\,u_3]\cap I}\e^{\mathrm{i}\phi(u)}\psi(u)\de u\right|\le
C\frac{\psi(\sigma+\sigma^\gamma/2)}{|\phi'(\sigma+\sigma^\gamma-\sigma^{\gamma/2})|}
\le C\frac{\sigma+\sigma^\gamma}{\sqrt{\sigma^{\gamma}/2}\sigma^{\gamma/2}}
\le C.
\]
As for $[u_3,\,u_4]$, we estimate it using the size of $\psi(u)$:
\begin{eqnarray*}
\left|\int_{[u_3,\,u_4]\cap I}\e^{\mathrm{i}\phi(u)}\psi(u)\de u\right|
&\le&\int_{u_3}^{\sigma+\sigma^\gamma+\sigma^{\gamma/2}}\psi(u)\de u
\le 2\sigma^{\gamma/2}\psi(u_3)\\
&\le& 2\sigma^{\gamma/2}\psi(\sigma+\sigma^{\gamma}/2)
\le\frac{2\sigma^{\gamma/2}}{\sqrt{\sigma^\gamma/2}}\le C.
\end{eqnarray*}
Finally, using Van der Corput's lemma again,
\[
\left|\int_{[u_4,\,\infty]\cap I}\e^{\mathrm{i}\phi(u)}\psi(u)\de u\right|\le
C\frac{\psi(\sigma+\sigma^\gamma+\sigma^{\gamma/2})}
{|\phi'(\sigma+\sigma^\gamma+\sigma^{\gamma/2})|}
\le C\frac{\sigma+\sigma^\gamma}{\sqrt{\sigma^\gamma+\sigma^{\gamma/2}}\sigma^{\gamma/2}}
\le C.
\]
This concludes the case $\sigma\ge 1$.
As for the remaining case, $0<\sigma\le 1$, we impose that the zero of $\phi'$
is $\sigma+\sigma^\gamma$, with $\gamma\le 1$ (when $\gamma$ grows from $-\infty$ to
$1$, $\sigma^\gamma$ decreases from $\infty$ to $\sigma$). Just as before, we have the
following estimates for $\phi'$ and $\psi$
\begin{eqnarray*}
\psi(u)&\le&\frac 1{\sqrt{u-\sigma}},\\
|\phi'(u)|&\ge&\frac{|\sigma+\sigma^\gamma-u|}{\sigma+\sigma^\gamma}
\ge\frac{|\sigma+\sigma^\gamma-u|}{2\sigma^\gamma}.
\end{eqnarray*}
Suppose $0\le\gamma\le 1$. Then
\[
\left|\int_{[2\sigma,\,3]\cap I}\e^{\mathrm{i}\phi(u)}\psi(u)\de u\right|\le
\int_{2\sigma}^{3}\psi(u)\de u\le C,
\]
and by Van der Corput's lemma,
\[
\left|\int_{[3,\,\infty]\cap I}\e^{\mathrm{i}\phi(u)}\psi(u)\de u\right|\le
C\frac{\psi(3)}{|\phi'(3)|}
\le C\frac{2\sigma^\gamma}{\sqrt{2}(3-\sigma-\sigma^{\gamma})}
\le C\sigma^\gamma\le C.
\]
If instead $\gamma<0$, then we divide the interval
$[2\sigma,\,\infty)$ into five subintervals, given by the
following partition
\begin{eqnarray*}
u_1&=&2\sigma,\\
u_2&=&3,\\
u_3&=&\max(3,\,\sigma+\sigma^\gamma/2),\\
u_4&=&\max(u_3,\,\sigma+\sigma^\gamma-\sigma^{\gamma/2}),\\
u_5&=&\sigma+\sigma^\gamma+\sigma^{\gamma/2},\\
\end{eqnarray*}
and study each case separately: the integrals along  the intervals
$[u_1,\,u_2]$ and $[u_4,\,u_5]$,
can be estimated by taking absolute values inside;
for the other intervals, apply Van der Corput's  lemma as usual.
\end{section}

\end{document}